\documentclass{gtart}

\def\ifplaintex{\expandafter\ifx\csname documentclass\endcsname\relax}

\def\gtp{{\mathsurround=0pt\it $\cal G\mskip-2mu$eometry \&\ 
$\cal T\!\!$opology $\cal P\!$ublications}}  % GT publications

\def\Addressesr{\bigskip
{\small \parskip 0pt \leftskip 0pt \rightskip 0pt plus 1fil \def\\{\par}
\sl\theaddress\par
\medskip
\rm Email:\stdspace\tt\theemail\hfill\rm Received:\qua\receiveddate \par}}

\def\recd{{\small Received:\qua\receiveddate\ifx\reviseddate\relax
\else\qquad Revised:\qua\reviseddate\fi\par}} 

\def\lognumber#1{\def\thelognumber{#1}}
\def\volumenumber#1{\def\thevolumenumber{#1}}
\def\volumeyear#1{\def\thevolumeyear{#1}}
\def\papernumber#1{\def\thepapernumber{#1}}
\def\pagenumbers#1#2{\def\startpage{#1}\def\finishpage{#2}}
\def\published#1{\def\publishdate{#1}}

\def\received#1{\def\receiveddate{#1}}

\def\accepted#1{\def\accepteddate{#1}}

\long\def\asciiabstract#1{\long\def\theasciiabstract{#1}}
\def\asciikeywords#1{\def\theasciikeywords{#1}}

\let\\\par\let\thelognumber\relax\let\thevolumenumber\relax
\let\thepapernumber\relax\let\thevolumeyear\relax\let\startpage\relax
\let\finishpage\relax\let\publishdate\relax\let\receiveddate\relax
\let\reviseddate\relax\let\accepteddate\relax\let\theasciititle\relax
\let\theasciiauthors\relax
\let\theasciiabstract\relax\let\theasciikeywords\relax

\let\theasciiemail\relax

\ifplaintex
\font\logobig=cmssbx10 scaled 3836
\font\logomed=cmssbx10 scaled 2557
\else
\font\logobig=cmssbx10 scaled 4200
\font\logomed=cmssbx10 scaled 2800
\fi

\long\def\makeagttitle{   %%% start of definition of \makeagttitle
\count0=\startpage
\agt\hfill      %   Journal title (top left) 
\hbox to 45truept{\vbox to 0pt{\vglue -13truept{\logomed A\kern -.37em{\logobig 
T}\kern -.38em G}\vss}\hss}
\break
{\small Volume \thevolumenumber\ (\thevolumeyear)
\startpage--\finishpage\nl
Published: \publishdate}

\vglue .25truein

{\parskip=0pt\leftskip 0pt plus
1fil\def\\{\par\smallskip}{\Large\bf\thetitle}\par\medskip} \vglue
0.05truein

{\parskip=0pt\leftskip 0pt plus 1fil\def\\{\par}{\sc\theauthors}
\par\medskip}%
 
\vglue 0.03truein 

{\small\leftskip 25truept\rightskip 25truept{\bf Abstract}\stdspace\theabstract

{\bf AMS Classification}\stdspace\theprimaryclass
\ifx\thesecondaryclass\relax\else; \thesecondaryclass\fi\par
{\bf Keywords}\stdspace \thekeywords\par}\vglue 7truept

}   %%%% end of definition of \makeagttitle

\ifplaintex
\font\phead=cmsl9 scaled 950
\font\pnum=cmbx10 scaled 913
\font\pfoot=cmsl9 scaled 950
\headline{\vbox to 0pt{\vskip -4.5mm\line{\small\phead\ifnum
\count0=\startpage ISSN 1472-2739 (on-line) 1472-2747 (printed)
\hfill {\pnum\folio}\else\ifodd\count0\def\\{ }% 
\ifx\theshorttitle\relax\thetitle\else\theshorttitle\fi\hfill{\pnum\folio}
\else\def\\{ and }{\pnum\folio}\hfill\ifx\theshortauthors\relax\theauthors
\else\theshortauthors\fi\fi\fi}\vss}}
\footline{\vbox to 0pt{\vglue 0mm\line{\small\pfoot\ifnum\count0=\startpage
\copyright\ \gtp\hfill\else
\agt, Volume \thevolumenumber\ (\thevolumeyear)\hfill\fi}\vss}}
\else
\font=cmsl9 scaled 1050
\font\lnum=cmbx10 
\font=cmsl9 scaled 1050
\makeatletter
\def\@oddhead{{\small\lhead\ifnum\count0=\startpage ISSN 1472-2739 
(on-line) 1472-2747 (printed)\hfill {\lnum\number\count0}\else\ifodd\count0
\def\\{ }\ifx\theshorttitle\relax \thetitle \else\theshorttitle\fi\hfill
{\lnum\number\count0}\else\def\\{ and }{\lnum\number\count0}
\hfill\ifx\theshortauthors\relax 
\theauthors\else\theshortauthors\fi\fi\fi}}\def\@evenhead{\@oddhead}
\def\@oddfoot{\small\lfoot\ifnum\count0=\startpage\copyright\ \gtp\hfill\else
\agt, Volume \thevolumenumber\ (\thevolumeyear)\hfill\fi}
\def\@evenfoot{\@oddfoot}
\makeatother
\fi
\let\maketitlepage\makeagttitle
\let\makeshorttitle\maketitlepage
\let\maketitle\maketitlepage

\newwrite\gtoutfile
\long\gdef\makeheadfile{  %%% start of definition of \makeheadfile
{\def\\{, }\def\s{ }
\immediate\openout\gtoutfile head.xxx
\immediate\write\gtoutfile{To: math@arxiv.org}
\immediate\write\gtoutfile{Subject: put OR rep NNNNN:ppppp}
\immediate\write\gtoutfile{--text follows this line--}
\immediate\write\gtoutfile{Proxy-for: \ifx\theasciiauthors\relax
\theauthors\else\theasciiauthors\fi\s<\ifx\theasciiemail\relax\theemail\else\theasciiemail\fi>}
\immediate\write\gtoutfile{\noexpand\\}
\immediate\write\gtoutfile{Authors: \ifx\theasciiauthors\relax
\theauthors\else\theasciiauthors\fi}
{\def\\{ }\immediate\write\gtoutfile{Title: \ifx\theasciititle\relax
\thetitle\else\theasciititle\fi}}
\immediate\write\gtoutfile{Subj-class: GT or SG, GR etc}
\immediate\write\gtoutfile{MSC-class: \theprimaryclass\ifx\thesecondaryclass\relax\else, \thesecondaryclass\fi}
\immediate\write\gtoutfile{Journal-ref: Algebr. Geom. Topol. \thevolumenumber\s
(\thevolumeyear) \startpage-\finishpage}
\immediate\write\gtoutfile{Comments: Published by Algebraic and
Geometric Topology at}
\immediate\write\gtoutfile{\s\s\s  http://www.maths.warwick.ac.uk/agt/AGTVol\thevolumenumber/agt-\thevolumenumber-\thepapernumber.abs.html}
\immediate\write\gtoutfile{\noexpand\\}
\immediate\write\gtoutfile{}
\ifx\theasciiabstract\relax
\immediate\write\gtoutfile{\theabstract}\else
\immediate\write\gtoutfile{\theasciiabstract}\fi
\immediate\write\gtoutfile{}
\immediate\write\gtoutfile{\noexpand\\}
\immediate\write\gtoutfile{}
\immediate\closeout\gtoutfile}}  %%% end of definition of \makeheadfile

\def\maketitlepage{\makeagttitle\makeheadfile}
\let\makeshorttitle\maketitlepage
\let\maketitle\maketitlepage

\def\ifplaintex{\expandafter\ifx\csname documentclass\endcsname\relax}

\def\gtp{{\mathsurround=0pt\it $\cal G\mskip-2mu$eometry \&\ 
$\cal T\!\!$opology $\cal P\!$ublications}}  % GT publications

\def\Addressesr{\bigskip
{\small \parskip 0pt \leftskip 0pt \rightskip 0pt plus 1fil \def\\{\par}
\sl\theaddress\par
\medskip
\rm Email:\stdspace\tt\theemail\hfill\rm Received:\qua\receiveddate \par}}

\def\recd{{\small Received:\qua\receiveddate\ifx\reviseddate\relax
\else\qquad Revised:\qua\reviseddate\fi\par}} 

\def\lognumber#1{\def\thelognumber{#1}}
\def\volumenumber#1{\def\thevolumenumber{#1}}
\def\volumeyear#1{\def\thevolumeyear{#1}}
\def\papernumber#1{\def\thepapernumber{#1}}
\def\pagenumbers#1#2{\def\startpage{#1}\def\finishpage{#2}}
\def\published#1{\def\publishdate{#1}}

\def\received#1{\def\receiveddate{#1}}

\def\accepted#1{\def\accepteddate{#1}}

\long\def\asciiabstract#1{\long\def\theasciiabstract{#1}}
\def\asciikeywords#1{\def\theasciikeywords{#1}}

\let\\\par\let\thelognumber\relax\let\thevolumenumber\relax
\let\thepapernumber\relax\let\thevolumeyear\relax\let\startpage\relax
\let\finishpage\relax\let\publishdate\relax\let\receiveddate\relax
\let\reviseddate\relax\let\accepteddate\relax\let\theasciititle\relax
\let\theasciiauthors\relax
\let\theasciiabstract\relax\let\theasciikeywords\relax

\let\theasciiemail\relax

\ifplaintex
\font\logobig=cmssbx10 scaled 3836
\font\logomed=cmssbx10 scaled 2557
\else
\font\logobig=cmssbx10 scaled 4200
\font\logomed=cmssbx10 scaled 2800
\fi

\long\def\makeagttitle{   %%% start of definition of \makeagttitle
\count0=\startpage
\agt\hfill      %   Journal title (top left) 
\hbox to 45truept{\vbox to 0pt{\vglue -13truept{\logomed A\kern -.37em{\logobig 
T}\kern -.38em G}\vss}\hss}
\break
{\small Volume \thevolumenumber\ (\thevolumeyear)
\startpage--\finishpage\nl
Published: \publishdate}

\vglue .25truein

{\parskip=0pt\leftskip 0pt plus
1fil\def\\{\par\smallskip}{\Large\bf\thetitle}\par\medskip} \vglue
0.05truein

{\parskip=0pt\leftskip 0pt plus 1fil\def\\{\par}{\sc\theauthors}
\par\medskip}%
 
\vglue 0.03truein 

{\small\leftskip 25truept\rightskip 25truept{\bf Abstract}\stdspace\theabstract

{\bf AMS Classification}\stdspace\theprimaryclass
\ifx\thesecondaryclass\relax\else; \thesecondaryclass\fi\par
{\bf Keywords}\stdspace \thekeywords\par}\vglue 7truept

}   %%%% end of definition of \makeagttitle

\ifplaintex
\font\phead=cmsl9 scaled 950
\font\pnum=cmbx10 scaled 913
\font\pfoot=cmsl9 scaled 950
\headline{\vbox to 0pt{\vskip -4.5mm\line{\small\phead\ifnum
\count0=\startpage ISSN 1472-2739 (on-line) 1472-2747 (printed)
\hfill {\pnum\folio}\else\ifodd\count0\def\\{ }% 
\ifx\theshorttitle\relax\thetitle\else\theshorttitle\fi\hfill{\pnum\folio}
\else\def\\{ and }{\pnum\folio}\hfill\ifx\theshortauthors\relax\theauthors
\else\theshortauthors\fi\fi\fi}\vss}}
\footline{\vbox to 0pt{\vglue 0mm\line{\small\pfoot\ifnum\count0=\startpage
\copyright\ \gtp\hfill\else
\agt, Volume \thevolumenumber\ (\thevolumeyear)\hfill\fi}\vss}}
\else
\font=cmsl9 scaled 1050
\font\lnum=cmbx10 
\font=cmsl9 scaled 1050
\makeatletter
\def\@oddhead{{\small\lhead\ifnum\count0=\startpage ISSN 1472-2739 
(on-line) 1472-2747 (printed)\hfill {\lnum\number\count0}\else\ifodd\count0
\def\\{ }\ifx\theshorttitle\relax \thetitle \else\theshorttitle\fi\hfill
{\lnum\number\count0}\else\def\\{ and }{\lnum\number\count0}
\hfill\ifx\theshortauthors\relax 
\theauthors\else\theshortauthors\fi\fi\fi}}\def\@evenhead{\@oddhead}
\def\@oddfoot{\small\lfoot\ifnum\count0=\startpage\copyright\ \gtp\hfill\else
\agt, Volume \thevolumenumber\ (\thevolumeyear)\hfill\fi}
\def\@evenfoot{\@oddfoot}
\makeatother
\fi
\let\maketitlepage\makeagttitle
\let\makeshorttitle\maketitlepage
\let\maketitle\maketitlepage

\newwrite\gtoutfile
\long\gdef\makeheadfile{  %%% start of definition of \makeheadfile
{\def\\{, }\def\s{ }
\immediate\openout\gtoutfile head.xxx
\immediate\write\gtoutfile{To: math@arxiv.org}
\immediate\write\gtoutfile{Subject: put OR rep NNNNN:ppppp}
\immediate\write\gtoutfile{--text follows this line--}
\immediate\write\gtoutfile{Proxy-for: \ifx\theasciiauthors\relax
\theauthors\else\theasciiauthors\fi\s<\ifx\theasciiemail\relax\theemail\else\theasciiemail\fi>}
\immediate\write\gtoutfile{\noexpand\\}
\immediate\write\gtoutfile{Authors: \ifx\theasciiauthors\relax
\theauthors\else\theasciiauthors\fi}
{\def\\{ }\immediate\write\gtoutfile{Title: \ifx\theasciititle\relax
\thetitle\else\theasciititle\fi}}
\immediate\write\gtoutfile{Subj-class: GT or SG, GR etc}
\immediate\write\gtoutfile{MSC-class: \theprimaryclass\ifx\thesecondaryclass\relax\else, \thesecondaryclass\fi}
\immediate\write\gtoutfile{Journal-ref: Algebr. Geom. Topol. \thevolumenumber\s
(\thevolumeyear) \startpage-\finishpage}
\immediate\write\gtoutfile{Comments: Published by Algebraic and
Geometric Topology at}
\immediate\write\gtoutfile{\s\s\s  http://www.maths.warwick.ac.uk/agt/AGTVol\thevolumenumber/agt-\thevolumenumber-\thepapernumber.abs.html}
\immediate\write\gtoutfile{\noexpand\\}
\immediate\write\gtoutfile{}
\ifx\theasciiabstract\relax
\immediate\write\gtoutfile{\theabstract}\else
\immediate\write\gtoutfile{\theasciiabstract}\fi
\immediate\write\gtoutfile{}
\immediate\write\gtoutfile{\noexpand\\}
\immediate\write\gtoutfile{}
\immediate\closeout\gtoutfile}}  %%% end of definition of \makeheadfile

\def\maketitlepage{\makeagttitle\makeheadfile}
\let\makeshorttitle\maketitlepage
\let\maketitle\maketitlepage

\def\ifplaintex{\expandafter\ifx\csname documentclass\endcsname\relax}

\def\gtp{{\mathsurround=0pt\it $\cal G\mskip-2mu$eometry \&\ 
$\cal T\!\!$opology $\cal P\!$ublications}}  % GT publications

\def\Addressesr{\bigskip
{\small \parskip 0pt \leftskip 0pt \rightskip 0pt plus 1fil \def\\{\par}
\sl\theaddress\par
\medskip
\rm Email:\stdspace\tt\theemail\hfill\rm Received:\qua\receiveddate \par}}

\def\recd{{\small Received:\qua\receiveddate\ifx\reviseddate\relax
\else\qquad Revised:\qua\reviseddate\fi\par}} 

\def\lognumber#1{\def\thelognumber{#1}}
\def\volumenumber#1{\def\thevolumenumber{#1}}
\def\volumeyear#1{\def\thevolumeyear{#1}}
\def\papernumber#1{\def\thepapernumber{#1}}
\def\pagenumbers#1#2{\def\startpage{#1}\def\finishpage{#2}}
\def\published#1{\def\publishdate{#1}}

\def\received#1{\def\receiveddate{#1}}

\def\accepted#1{\def\accepteddate{#1}}

\long\def\asciiabstract#1{\long\def\theasciiabstract{#1}}
\def\asciikeywords#1{\def\theasciikeywords{#1}}

\let\\\par\let\thelognumber\relax\let\thevolumenumber\relax
\let\thepapernumber\relax\let\thevolumeyear\relax\let\startpage\relax
\let\finishpage\relax\let\publishdate\relax\let\receiveddate\relax
\let\reviseddate\relax\let\accepteddate\relax\let\theasciititle\relax
\let\theasciiauthors\relax
\let\theasciiabstract\relax\let\theasciikeywords\relax

\let\theasciiemail\relax

\ifplaintex
\font\logobig=cmssbx10 scaled 3836
\font\logomed=cmssbx10 scaled 2557
\else
\font\logobig=cmssbx10 scaled 4200
\font\logomed=cmssbx10 scaled 2800
\fi

\long\def\makeagttitle{   %%% start of definition of \makeagttitle
\count0=\startpage
\agt\hfill      %   Journal title (top left) 
\hbox to 45truept{\vbox to 0pt{\vglue -13truept{\logomed A\kern -.37em{\logobig 
T}\kern -.38em G}\vss}\hss}
\break
{\small Volume \thevolumenumber\ (\thevolumeyear)
\startpage--\finishpage\nl
Published: \publishdate}

\vglue .25truein

{\parskip=0pt\leftskip 0pt plus
1fil\def\\{\par\smallskip}{\Large\bf\thetitle}\par\medskip} \vglue
0.05truein

{\parskip=0pt\leftskip 0pt plus 1fil\def\\{\par}{\sc\theauthors}
\par\medskip}%
 
\vglue 0.03truein 

{\small\leftskip 25truept\rightskip 25truept{\bf Abstract}\stdspace\theabstract

{\bf AMS Classification}\stdspace\theprimaryclass
\ifx\thesecondaryclass\relax\else; \thesecondaryclass\fi\par
{\bf Keywords}\stdspace \thekeywords\par}\vglue 7truept

}   %%%% end of definition of \makeagttitle

\ifplaintex
\font\phead=cmsl9 scaled 950
\font\pnum=cmbx10 scaled 913
\font\pfoot=cmsl9 scaled 950
\headline{\vbox to 0pt{\vskip -4.5mm\line{\small\phead\ifnum
\count0=\startpage ISSN 1472-2739 (on-line) 1472-2747 (printed)
\hfill {\pnum\folio}\else\ifodd\count0\def\\{ }% 
\ifx\theshorttitle\relax\thetitle\else\theshorttitle\fi\hfill{\pnum\folio}
\else\def\\{ and }{\pnum\folio}\hfill\ifx\theshortauthors\relax\theauthors
\else\theshortauthors\fi\fi\fi}\vss}}
\footline{\vbox to 0pt{\vglue 0mm\line{\small\pfoot\ifnum\count0=\startpage
\copyright\ \gtp\hfill\else
\agt, Volume \thevolumenumber\ (\thevolumeyear)\hfill\fi}\vss}}
\else
\font=cmsl9 scaled 1050
\font\lnum=cmbx10 
\font=cmsl9 scaled 1050
\makeatletter
\def\@oddhead{{\small\lhead\ifnum\count0=\startpage ISSN 1472-2739 
(on-line) 1472-2747 (printed)\hfill {\lnum\number\count0}\else\ifodd\count0
\def\\{ }\ifx\theshorttitle\relax \thetitle \else\theshorttitle\fi\hfill
{\lnum\number\count0}\else\def\\{ and }{\lnum\number\count0}
\hfill\ifx\theshortauthors\relax 
\theauthors\else\theshortauthors\fi\fi\fi}}\def\@evenhead{\@oddhead}
\def\@oddfoot{\small\lfoot\ifnum\count0=\startpage\copyright\ \gtp\hfill\else
\agt, Volume \thevolumenumber\ (\thevolumeyear)\hfill\fi}
\def\@evenfoot{\@oddfoot}
\makeatother
\fi
\let\maketitlepage\makeagttitle
\let\makeshorttitle\maketitlepage
\let\maketitle\maketitlepage

\newwrite\gtoutfile
\long\gdef\makeheadfile{  %%% start of definition of \makeheadfile
{\def\\{, }\def\s{ }
\immediate\openout\gtoutfile head.xxx
\immediate\write\gtoutfile{To: math@arxiv.org}
\immediate\write\gtoutfile{Subject: put OR rep NNNNN:ppppp}
\immediate\write\gtoutfile{--text follows this line--}
\immediate\write\gtoutfile{Proxy-for: \ifx\theasciiauthors\relax
\theauthors\else\theasciiauthors\fi\s<\ifx\theasciiemail\relax\theemail\else\theasciiemail\fi>}
\immediate\write\gtoutfile{\noexpand\\}
\immediate\write\gtoutfile{Authors: \ifx\theasciiauthors\relax
\theauthors\else\theasciiauthors\fi}
{\def\\{ }\immediate\write\gtoutfile{Title: \ifx\theasciititle\relax
\thetitle\else\theasciititle\fi}}
\immediate\write\gtoutfile{Subj-class: GT or SG, GR etc}
\immediate\write\gtoutfile{MSC-class: \theprimaryclass\ifx\thesecondaryclass\relax\else, \thesecondaryclass\fi}
\immediate\write\gtoutfile{Journal-ref: Algebr. Geom. Topol. \thevolumenumber\s
(\thevolumeyear) \startpage-\finishpage}
\immediate\write\gtoutfile{Comments: Published by Algebraic and
Geometric Topology at}
\immediate\write\gtoutfile{\s\s\s  http://www.maths.warwick.ac.uk/agt/AGTVol\thevolumenumber/agt-\thevolumenumber-\thepapernumber.abs.html}
\immediate\write\gtoutfile{\noexpand\\}
\immediate\write\gtoutfile{}
\ifx\theasciiabstract\relax
\immediate\write\gtoutfile{\theabstract}\else
\immediate\write\gtoutfile{\theasciiabstract}\fi
\immediate\write\gtoutfile{}
\immediate\write\gtoutfile{\noexpand\\}
\immediate\write\gtoutfile{}
\immediate\closeout\gtoutfile}}  %%% end of definition of \makeheadfile

\def\maketitlepage{\makeagttitle\makeheadfile}
\let\makeshorttitle\maketitlepage
\let\maketitle\maketitlepage

\def\ifplaintex{\expandafter\ifx\csname documentclass\endcsname\relax}

\def\gtp{{\mathsurround=0pt\it $\cal G\mskip-2mu$eometry \&\ 
$\cal T\!\!$opology $\cal P\!$ublications}}  % GT publications

\def\Addressesr{\bigskip
{\small \parskip 0pt \leftskip 0pt \rightskip 0pt plus 1fil \def\\{\par}
\sl\theaddress\par
\medskip
\rm Email:\stdspace\tt\theemail\hfill\rm Received:\qua\receiveddate \par}}

\def\recd{{\small Received:\qua\receiveddate\ifx\reviseddate\relax
\else\qquad Revised:\qua\reviseddate\fi\par}} 

\def\lognumber#1{\def\thelognumber{#1}}
\def\volumenumber#1{\def\thevolumenumber{#1}}
\def\volumeyear#1{\def\thevolumeyear{#1}}
\def\papernumber#1{\def\thepapernumber{#1}}
\def\pagenumbers#1#2{\def\startpage{#1}\def\finishpage{#2}}
\def\published#1{\def\publishdate{#1}}

\def\received#1{\def\receiveddate{#1}}

\def\accepted#1{\def\accepteddate{#1}}

\long\def\asciiabstract#1{\long\def\theasciiabstract{#1}}
\def\asciikeywords#1{\def\theasciikeywords{#1}}

\let\\\par\let\thelognumber\relax\let\thevolumenumber\relax
\let\thepapernumber\relax\let\thevolumeyear\relax\let\startpage\relax
\let\finishpage\relax\let\publishdate\relax\let\receiveddate\relax
\let\reviseddate\relax\let\accepteddate\relax\let\theasciititle\relax
\let\theasciiauthors\relax
\let\theasciiabstract\relax\let\theasciikeywords\relax

\let\theasciiemail\relax

\ifplaintex
\font\logobig=cmssbx10 scaled 3836
\font\logomed=cmssbx10 scaled 2557
\else
\font\logobig=cmssbx10 scaled 4200
\font\logomed=cmssbx10 scaled 2800
\fi

\long\def\makeagttitle{   %%% start of definition of \makeagttitle
\count0=\startpage
\agt\hfill      %   Journal title (top left) 
\hbox to 45truept{\vbox to 0pt{\vglue -13truept{\logomed A\kern -.37em{\logobig 
T}\kern -.38em G}\vss}\hss}
\break
{\small Volume \thevolumenumber\ (\thevolumeyear)
\startpage--\finishpage\nl
Published: \publishdate}

\vglue .25truein

{\parskip=0pt\leftskip 0pt plus
1fil\def\\{\par\smallskip}{\Large\bf\thetitle}\par\medskip} \vglue
0.05truein

{\parskip=0pt\leftskip 0pt plus 1fil\def\\{\par}{\sc\theauthors}
\par\medskip}%
 
\vglue 0.03truein 

{\small\leftskip 25truept\rightskip 25truept{\bf Abstract}\stdspace\theabstract

{\bf AMS Classification}\stdspace\theprimaryclass
\ifx\thesecondaryclass\relax\else; \thesecondaryclass\fi\par
{\bf Keywords}\stdspace \thekeywords\par}\vglue 7truept

}   %%%% end of definition of \makeagttitle

\ifplaintex
\font\phead=cmsl9 scaled 950
\font\pnum=cmbx10 scaled 913
\font\pfoot=cmsl9 scaled 950
\headline{\vbox to 0pt{\vskip -4.5mm\line{\small\phead\ifnum
\count0=\startpage ISSN 1472-2739 (on-line) 1472-2747 (printed)
\hfill {\pnum\folio}\else\ifodd\count0\def\\{ }% 
\ifx\theshorttitle\relax\thetitle\else\theshorttitle\fi\hfill{\pnum\folio}
\else\def\\{ and }{\pnum\folio}\hfill\ifx\theshortauthors\relax\theauthors
\else\theshortauthors\fi\fi\fi}\vss}}
\footline{\vbox to 0pt{\vglue 0mm\line{\small\pfoot\ifnum\count0=\startpage
\copyright\ \gtp\hfill\else
\agt, Volume \thevolumenumber\ (\thevolumeyear)\hfill\fi}\vss}}
\else
\font=cmsl9 scaled 1050
\font\lnum=cmbx10 
\font=cmsl9 scaled 1050
\makeatletter
\def\@oddhead{{\small\lhead\ifnum\count0=\startpage ISSN 1472-2739 
(on-line) 1472-2747 (printed)\hfill {\lnum\number\count0}\else\ifodd\count0
\def\\{ }\ifx\theshorttitle\relax \thetitle \else\theshorttitle\fi\hfill
{\lnum\number\count0}\else\def\\{ and }{\lnum\number\count0}
\hfill\ifx\theshortauthors\relax 
\theauthors\else\theshortauthors\fi\fi\fi}}\def\@evenhead{\@oddhead}
\def\@oddfoot{\small\lfoot\ifnum\count0=\startpage\copyright\ \gtp\hfill\else
\agt, Volume \thevolumenumber\ (\thevolumeyear)\hfill\fi}
\def\@evenfoot{\@oddfoot}
\makeatother
\fi
\let\maketitlepage\makeagttitle
\let\makeshorttitle\maketitlepage
\let\maketitle\maketitlepage

\newwrite\gtoutfile
\long\gdef\makeheadfile{  %%% start of definition of \makeheadfile
{\def\\{, }\def\s{ }
\immediate\openout\gtoutfile head.xxx
\immediate\write\gtoutfile{To: math@arxiv.org}
\immediate\write\gtoutfile{Subject: put OR rep NNNNN:ppppp}
\immediate\write\gtoutfile{--text follows this line--}
\immediate\write\gtoutfile{Proxy-for: \ifx\theasciiauthors\relax
\theauthors\else\theasciiauthors\fi\s<\ifx\theasciiemail\relax\theemail\else\theasciiemail\fi>}
\immediate\write\gtoutfile{\noexpand\\}
\immediate\write\gtoutfile{Authors: \ifx\theasciiauthors\relax
\theauthors\else\theasciiauthors\fi}
{\def\\{ }\immediate\write\gtoutfile{Title: \ifx\theasciititle\relax
\thetitle\else\theasciititle\fi}}
\immediate\write\gtoutfile{Subj-class: GT or SG, GR etc}
\immediate\write\gtoutfile{MSC-class: \theprimaryclass\ifx\thesecondaryclass\relax\else, \thesecondaryclass\fi}
\immediate\write\gtoutfile{Journal-ref: Algebr. Geom. Topol. \thevolumenumber\s
(\thevolumeyear) \startpage-\finishpage}
\immediate\write\gtoutfile{Comments: Published by Algebraic and
Geometric Topology at}
\immediate\write\gtoutfile{\s\s\s  http://www.maths.warwick.ac.uk/agt/AGTVol\thevolumenumber/agt-\thevolumenumber-\thepapernumber.abs.html}
\immediate\write\gtoutfile{\noexpand\\}
\immediate\write\gtoutfile{}
\ifx\theasciiabstract\relax
\immediate\write\gtoutfile{\theabstract}\else
\immediate\write\gtoutfile{\theasciiabstract}\fi
\immediate\write\gtoutfile{}
\immediate\write\gtoutfile{\noexpand\\}
\immediate\write\gtoutfile{}
\immediate\closeout\gtoutfile}}  %%% end of definition of \makeheadfile

\def\maketitlepage{\makeagttitle\makeheadfile}
\let\makeshorttitle\maketitlepage
\let\maketitle\maketitlepage

\lognumber{15}
\volumenumber{2}
\volumeyear{2002}
\papernumber{15}
\published{10 May 2002}
\pagenumbers{311}{336}
\received{18 February 2002}
\accepted{8 May 2002}

\usepackage{amsmath,amscd,amssymb,epsf,pb-diagram,lamsarrow,pb-lams}

\newtheorem{thm}{Theorem}[section]  

\newtheorem{lem}[thm]{Lemma}

\newtheorem{add}[thm]{Addendum}
\newtheorem*{claim}{Claim}

\newtheorem{bigthm}{Theorem}
\newtheorem{bigcor}[bigthm]{Corollary}

\theoremstyle{definition}
\newtheorem{defn}[thm]{Definition}

\theoremstyle{definition}

\theoremstyle{remark}
\newtheorem{rem}[thm]{Remark}

\newtheorem{exs}[thm]{Examples}

\def\Top{\bold T\bold o \bold p}
\def\vo{\varOmega}
\def\vs{\varSigma}
\def\smsh{\wedge}

\def\id{\text{id}}
\def\dbslash{/\!\! /}
\def\codim{\text{\rm codim\,}}
\def\:{\colon\thinspace}
\def\holim{\text{holim\,}}
\def\bold{\mathbf}
\def\Bbb{\mathbb}
\def\cal{\mathcal}

\begin{document}
\title{Embedding, compression and\\fiberwise homotopy theory}
\author{John R. Klein}
\address{Department of Mathematics, Wayne State University\\Detroit, 
MI 48202, USA}
\email{klein@math.wayne.edu}
\begin{abstract}Given Poincar\'e spaces $M$ and $X$,
we study the possibility of 
compressing embeddings of $M {\times} I$ in $X {\times}I$ 
down to embeddings of $M$ in $X$. This results in
a new approach to embedding in the metastable range
both in the smooth and Poincar\'e duality categories.
\end{abstract}
\asciiabstract{Given Poincare spaces M and X,
we study the possibility of 
compressing embeddings of M x I in X x I 
down to embeddings of M in X. This results in
a new approach to embedding in the metastable range
both in the smooth and Poincare duality categories.}
\asciikeywords{Poincare space, embedding, fiberwise homotopy}
\maketitle

\section{Introduction} 
Let $M$ and $X$ be compact $n$-manifolds.
The word compression of the title refers to a situation in
which one  is presented with  an embedding of
$M{\times} I$ in the interior of $ X {\times} I$
and  then tries
to decide whether it arises from an embedding 
of $M$ in  $X$, up to isotopy. If so, then the original embedding
{\it compresses.} One aim of the present paper is to decide
when this is possible.

The compression problem is mirrored in the Poincar\'e duality category.
From now on, let $M$ and $X$ be Poincar\'e duality spaces of dimension $n$. 
One says that $M$ {\it 
(Poincar\'e) embeds}
in $X$ with {\it complement} $C$ if there exists a 
decomposition $X \simeq  M \cup_{\partial M} C$
in which 
$\partial M\amalg \partial X$ is identified with a 
Poincar\'e duality boundary for
$C$ (we also assume a compatibility of fundamental classes---see 
\ref{defn-embed} below.) 

It will be convenient to have separate notation
for intervals of different lengths. 
Let $I = [0,1]$ and $J = [1/3,2/3]$. For a subspace  $S \subset  I$ set
$M_S := M {\times} S$.
We start with the following data: an
embedding of the $(n{+}1)$-dimensional Poincar\'e space 
$M_J$ in $X_I$ with
complement $W$. This gives us a map $\nu \:M
\to W$ by taking the composition 
$$
M_{1/3} \subset
\partial M_J \to W \, .
$$

Let $R(X)$ denote the category of {\it retractive spaces} over $X$.
An object of $R(X)$ is a space $Y$ equipped with
maps $s_Y\:X \to Y$ and $r_Y\:Y \to X$ 
(called respectively {\it inclusion} and {\it retraction})
such that $r_Y\circ s_Y$ is
the identity (objects are usually specified without reference
to their structure maps). 
A morphism $Y \to Z$ is a map of spaces which is compatible 
with the structure maps. According to Quillen \cite{Quillen}, $R(X)$ is
a model category in which a {\it weak equivalence} is a morphism
$Y \to Z$ which when considered as a map of
 spaces is a weak homotopy equivalence
(for the remaining structure, see \ref{fib-space} below). Hence,
it makes sense to speak of its {\it homotopy category} $\text{ho}R(X)$.

The inclusion  $X_0 \subset W$ and the
composite $W \to X_I\to X$ equip the space $W$ with
the structure of an object of $R(X)$. Let $M^+$ denote the
object of $R(X)$ given by taking the disjoint union of 
$M$ with $X$;
the inclusion $X \to M^+$ is evident
and the retraction $M^+ \to X$ is defined to be the composite
$$
\begin{CD}
M \amalg X = M_{1/3} \amalg X_0  \subset M_J \amalg X_0 \to
X_I @> \text{project} >> X\,\, .
\end{CD}
$$

With respect to these conventions, the map $\nu\:M \to W$ induces
a morphism 
$$
\nu^+\:M^+ \to W
$$ 
of $R(X)$.
Then $\nu^+$ determines a fiberwise
 homotopy class 
$$
[\nu^+] \in [M^+,W]_X \, . 
$$

\begin{rem}This will be the primary obstruction to compression. Informally,
it should be thought of as  measuring  the self-linking of $M$
in $X_I$. Several authors have studied non-fiberwise versions of this
construction (see Hirsch \cite{Hirsch}, Levitt
\cite{Levitt2} and  Williams \cite{Williams2}).
\end{rem}

Following Goodwillie \cite{Good3}, the  {\it homotopy codimension} of $M$
is $\ge q$, if 
\begin{itemize}
\item 
$M$ is homotopy equivalent to a CW complex of $\dim \le n{-}q$, and 
\item
the inclusion  $\partial M \to M$ is $(q{-}1)$-connected. 
\end{itemize}
In what follows, we write $\codim M \ge q$.
By a result of Wall \cite{Wall1},  the 
first condition is a consequence of the second whenever $q \ge 3$.

\begin{exs}
(1)\qua If $M$ is regular neighborhood of $p$-dimensional complex in an
$n$-dimensional manifold, then $\codim M \ge n{-}p$. 
 \smallskip

(2)\qua Let $V^p$ be a closed Poincar\'e space of
dimension $p$ equipped with  an $(n{-}p{-}1)$-spherical
fibration $\xi\: S(\xi) \to V$. 
Let $D(\xi)$ be the mapping cylinder of $\xi$. Then
$(D(\xi),S(\xi))$ is a Poincar\'e pair of dimension $n$ with
$\codim D(\xi) \ge n{-}p$.
\end{exs}

We now state the main result.

\begin{bigthm} Assume 
$\codim M \ge n{-}p \ge 3$ and  $3p{+}4 \le 2n$.
Then there exists an embedding of $M$ in
$X$ which induces the given embedding of $M_J$
in $X_I$ (up to ``concordance'') if and only if
$[\nu^+] \in [M^+,W]_X$ is trivial.
\end{bigthm}

We remark that this is valid in both the smooth
and Poincar\'e cases (the smooth case follows by application
of the surgery  machine---see below). 
In the special case $X = D^n$ is a disk, 
Theorem A reduces to a non-fiberwise result which
is implicit in the work of Williams \cite{Williams}. In fact,
our proof of Theorem A is a fiberwisation of 
one of Williams' arguments.

With respect to the numerical assumptions of Theorem A, we have

\begin{add}\label{stab} The map of fiberwise homotopy classes  
$$
\vs_X \: [M^+,W]_X \to [\vs_X M^+,\vs_X W]_X
$$
is an isomorphism, where $\vs_X$ denotes fiberwise suspension.
Consequently, the obstruction to compression $[\nu^+]$ is stable.
\end{add}

 This is proved in \S7
using  the Freudenthal suspension theorem
for $\text{ho}R(X)$ (cf.\ \ref{freud} below).

\subsection{Unstable fiberwise normal invariants}
Let $M$ and $X$ be $n$-dimensional
Poincar\'e spaces, and let $f\:M \to X$ be a map. 
These data define an object
$$
M\dbslash \partial M \in R(X)
$$ 
 whose underlying space
is $X \cup_{f|\partial M} M$ (note: collapsing $X$ to
 a point gives the quotient $M/\partial M$).
Similarly, we have $X\dbslash \partial X \in R(X)$ 
which turns out to be the {\it double}
$X\cup_{\partial X} X$ (which gives
$X^+$ if $\partial X$ is empty.)

If
$f\: M \to X$ is the underlying map of an embedding of $M$
in $X$ with complement $C$, then there is 
an associated fiberwise homotopy
class 
$$
\alpha_f \in [X\dbslash \partial X,M\dbslash \partial M]_X
$$
defined by taking 
$$
\begin{CD}
X \cup_{\partial X}  X @< \simeq << 
X \cup_{\partial X} (C \cup_{\partial M} M) @>>>
X \cup_X (X \cup_{\partial M} M) = M\dbslash \partial M\, .
\end{CD}
$$
This is the {\it fiberwise (Thom-Pontryagin) collapse} of the embedding.

By analogy
with Smale-Hirsch theory, a map
$f\:M \to X$ is said to
{\it (Poincar\'e) immerse} if there exists an integer $j\ge 0$ such that 
$f{\times}\id\: M {\times} D^j \to X {\times} D^j$ is the underlying
map of some embedding.

\begin{rem}  A fact we won't need, but which is
nevertheless true, is that $f$ Poincar\'e immerses
if and only if there is a stable fiber homotopy equivalence $f^*\nu_X \simeq \nu_M$, where $\nu_X$ and $\nu_M$ denote
the Spivak normal fibrations of $X$ and $M$ respectively.
For a proof of this, see \cite{Klein5}.
\end{rem}

Taking the fiberwise collapse of the embedding
$M{\times} D^j \to X {\times}D^j$ enables us to associate
a fiberwise {\it stable} homotopy
class
$$
\alpha^{\text{st}}_f \in \{ X\dbslash \partial X,M\dbslash \partial M\}_X
$$
called the {\it fiberwise (stable) normal invariant} of the 
immersion (this is independent of the choice of embedding.)

Obviously, a necessary obstruction to compressing 
the given embedding to  an embedding of $M$ in $X$ is
that $\alpha^{\text{st}}_f$ should  desuspend to an element $\alpha_f \in
[ X\dbslash \partial X,M\dbslash \partial M]_X$.
Call any such desuspension a {\it fiberwise unstable normal invariant}
of the immersion.

\begin{bigthm} Assume $f\:M \to X$  immerses. 
Again, suppose that $\codim M \ge n{-} p \ge 3$ 
and $3p{+}4 \le 2n$. Then $f$ 
embeds  (inducing the given immersion) if and
only if there exists a fiberwise unstable normal invariant
$\alpha_f$.
Moreover, the embedding can be chosen so that its collapse induces 
$\alpha_f$.
\end{bigthm}

In the case  $\partial X = \emptyset$, Richter has
also proved  Theorem B  using fiberwise Hopf
invariants and fiberwise $S$-duality. By contrast, we will deduce Theorem B from Theorem A (in fact, the  theorems are equivalent).

A consequence of the above is
a  Whitney embedding theorem for immersions in the
Poincar\'e duality category:

\begin{bigcor} Assume  $f\:M^p \to X^n$ immerses, where
$\codim M \ge n {-} p\ge 3$ and $2p {+} 1 \le n$. Then 
$f$ embeds (inducing the given immersion up to concordance).
\end{bigcor}

This follows because the fiberwise stable normal invariant 
destabilizes by \ref{freud}.

\subsection{A Levine style embedding theorem}
When $X$ is `highly' connected, Theorem B simplifies to a 
non-fiberwise statement. 
Here is its formulation: given an immersion of $f\:M \to X$ 
as above,
there is an associated stable (Thom-Pontryagin) collapse
$$
\beta^{\text{st}} \in \{ X/\partial X, M/\partial M\} \, .
$$
Any homotopy class
$$
\beta \in [X/\partial X,M/\partial M]
$$
which suspends to $\beta^{\text{st}}$ is called  
an {\it unstable normal invariant}. 

\begin{bigthm} Assume $\codim M \ge n{-}p \ge 3$,
$X$ is $[p/2]$-connected and $3p{+}4 \le 2n$. 
Then there exists an embedding inducing the given immersion of 
$M$ in $X$ if and only if there exists 
an unstable normal invariant $\beta$.
Moreover, the embedding can be chosen 
so that its collapse
coincides with $\beta$.
\end{bigthm}

For example, if we take $X = D^n$ then we recover the Williams-Richter
embedding theorem \cite{Williams}, \cite{Richter}. Levine's 
embedding theorem \cite[Thm.\ 4]{Levine} amounts to
the case when $X$ is a smooth $n$-manifold and $M = D(\xi)$
is the unit disk bundle of a vector bundle over a smooth
$p$-manifold $V$.

\subsection{Embedding  spheres in the middle dimension}
In applications to surgery on Poincar\'e spaces, 
one of the main issues is whether or not homotopy classes in
the middle dimension are represented by 
`framed' embedded spheres.

Let $X^n$ be a Poincar\'e space, and suppose that $n = 2p$.
Set 
$$
P := S^p {\times} D^p\, , $$
 and suppose that
$f\: P \to  X$ is a map which immerses. Let 
$\widetilde X$
be the universal cover of $X$, and let $\pi$ be the group of
deck transformations. A map $Y \to X$ then induces a $\pi$-covering
of $Y$ which we denote by $\widetilde Y$. Note that
$\widetilde X/\partial \widetilde X$ is a based $\pi$-space,
which is free in the based sense.

The immersion $f$ gives rise to an equivariant stable homotopy class
$$
\widetilde \beta^{\text{st}} \in \{ \widetilde X/\partial \widetilde X,
\widetilde P/\partial \widetilde P\}^\pi \, ,
$$
called the {\it equivariant stable collapse}.
This  is constructed as follows: choose a representative 
embedding for  $f{\times} \id_{D^j}\: P {\times} D^j \to X {\times} D^j$.
The diagram for this embedding can then be pulled-back along 
$\widetilde X$. The Thom-Pontryagin collapse of the resulting
diagram of $\pi$-spaces then yields
$\widetilde \beta^{\text{st}}$.

\begin{bigthm} \label{spheres} Assume $p > 2$. An immersion $f\: P \to X$
is represented by an embedding if and only if 
the equivariant stable collapse desuspends to an element
 $\widetilde \beta \in [\widetilde X/\partial \widetilde X,
\widetilde P/\partial \widetilde P]^\pi$. Furthermore, 
the embedding can be chosen so that its equivariant collapse is
$\widetilde \beta$.
\end{bigthm}

\subsection{Embedded thickenings}
Up until now, we have discussed embedding theorems
between Poincar\'e spaces having the same dimension. 
In a previous paper \cite{Klein}, we studied the following related problem:
suppose that $K$ is the homotopy type of a finite complex,
 $X^n$ is a Poincar\'e space, and $f\: K \to X$
is a map. Does there exist a `Poincar\'e boundary' for $K$, say
$A \to K$,  such that $f\: K \to X$ embeds? (More precisely,
we should really replace $K$ by the mapping cylinder of the
map $A \to K$ to get a Poincar\'e pair.)  Additionally,
one assumes a codimensional restriction: $k \le n{-}3$,
where $k$ is the {\it homotopy dimension} of  $K$ (an integer
such that $K$ is homotopy equivalent to a CW complex of that dimension).

This is the notion of Poincar\'e embedding in which the `normal data'
are not {\it a priori} chosen. In \cite{Klein} we termed this notion
a {\it PD embedding}. In this paper, we will
call it an {\it embedded thickening}, since
the choice of Poincar\'e boundary is a `Poincar\'e thickening'
of $K$.

An important special case of this concept
is when $K$ itself is a closed Poincar\'e space. In this instance,
the homotopy fiber of the map $A \to K$ is a sphere, and
one recovers the notion of Poincar\'e embedding used by Wall
\cite[Chap.\ 11]{Wall}. 

In \cite{Klein}, we proved that $f\: K^k \to X^n$ 
embedded thickens whenever $f$ is $(2k{-}$ $n{+}2)$-connected.
It was expected that this is not the sharpest result,
for in the smooth case, this result can be improved by one dimension.
We show that the result can be improved by one dimension in the
range of Theorem A:

\begin{bigthm}  Assume 
$f\: K \to X$ is $(2k{-}n{+}1)$-connected,
$k\le n{-}3$ and $3k{+}4 \le 2n$. Then there exists an embedded
thickening of $f$.
\end{bigthm}

Note that this immediately implies
the Poincar\'e versions  of the `easy' and `hard'
Whitney embedding theorems:
let $f\:K \to X$ be a map with $k \le n{-}3$.

\begin{bigcor}\label{hard-whitney}{\rm(1)}\qua If $2k{+}1 \le n$, then 
$f$ embedded thickens.
\smallskip

{\rm(2)}\qua  If $2k \le n$ and  
and $f$ is $1$-connected, then 
$f$ embedded thickens with the  possible exception of the 
case  $k {=} 3$ and  $n {=} 6$.
\end{bigcor}

\begin{rem} The first part of the corollary settles an issue raised by
 Levitt \cite[p.\ 402]{Levitt2}.
\end{rem}

Another application yields 
an extension of \cite[Cor.\ C]{Klein}, which concerns the existence
of the unstable homotopy tangent bundle for Poincar\'e spaces:

\begin{bigcor} Let $X^n$ be a $1$-connected closed Poincar\'e space. Then
the diagonal $X \to X{\times} X$ has an embedded thickening.
\end{bigcor}

This follows by Theorem F  if $n \ge 4$, and is trivial if
$n <  4$.

\subsection{Smooth embeddings}
If $M$ and $X$ are compact smooth manifolds, then
the Browder-Casson-Sullivan-Wall theorem \cite[Chap.\ 11]{Wall}
shows that
all of the above results imply smooth embedding results,
(some new, some known). We leave it to the reader to make
sense of this translation.

The inequality $3p {+} 4 \le 2n$ can be improved to $3p {+} 3 \le 2n$ 
in the smooth case: in proving Theorem A we make use of the 
relative embedding theorem of \cite{Klein3}, 
which is the Poincar\'e variant of a result of Hodgson
\cite{Hodgson2} 
with a loss of one dimension. In the smooth case, 
Hodgson's result may be directly
substituted in the appropriate part of the proof
of Theorem A  to yield the sharper result.

\subsection{History}
The concept of Poincar\'e embedding surfaced in an
attempt to understand smooth embeddings within the framework
of surgery theory. The Browder-Casson-Sullivan-Wall theorem asserts that the 
smooth embedding problem of $M^n$ in $X^n$  is equivalent to the corresponding
Poincar\'e embedding problem as long as $n \ge 6$ and $\codim M \ge 3$.
Consequently, the problem of smooth embedding is reduced to
homotopy theory.

The inequality  $3p {+}3 \le 2n$ is called the {\it metastable} range. Roughly,
it is the place where triple point obstructions don't arise for dimensional
reasons. 

From 1960-1975 there emerged (at least) 
three different strategies to
(smooth) embedding in the metastable range.
Firstly, there was the school of Haefliger, which reduced the problem
to a question about isovariant maps $M^{{\times} 2} \to X^{{\times} 2}$
(an equivariant 
 map such that the inverse image of the diagonal of $X$ coincides
with the diagonal of $M$). Secondly, there was the bordism theoretic
approach, as seen in the papers of Dax \cite{Dax} and 
Hatcher-Quinn \cite{Hatch-Quinn}. Both of these schools relied heavily on
the Whitney trick and/or engulfing methods.

Lastly, there was
the surgery school---most notably the works of 
Browder \cite{Browder2}, \cite{Browder3} and  Wall 
\cite{Wall}---which reduced the problem of smooth
embedding to that of Poincar\'e embedding. This approach began with
Levine \cite{Levine}, who, using  surgery, 
constructed embeddings
from unstable normal invariants 
when the source $M$ is the total space of a  disk bundle over a
smooth manifold and
the ambient space $X$ is an $n$-sphere, or more generally when
$X$ is a sufficiently highly connected manifold. Here, the role of the normal
bundle is prominent.

Later,  Williams \cite{Williams}, \cite{Williams2},
Rigdon-Williams \cite{Rigdon-Williams}  
and Richter \cite{Richter}, extended Levine's work to
the case when $M$ is a Poincar\'e space and $X = D^n$. The work
of Williams \it et.\ al.\ \rm used smooth manifold techniques to
deduce results about Poincar\'e embeddings. Richter gave the first
manifold-free proof of Williams' results using homotopy theory.

It was only recently observed \cite{Klein} that
 fiberwise homotopy theory technology
 was to play a role in extending the 
surgery approach to an arbitrary ambient Poincar\'e space $X$. 
This connection was discovered  by 
Shmuel Weinberger and  the author (independently).  The present work
is an attempt to complete  the thread begun by the surgery school.  
 
\subsection{Outline}  Section 2 is mostly language; the reader should 
be familiar with the majority of material in this section. 
In \S3 we show that the existence of a fiberwise
normal invariant is sufficient to give
an embedding of $M_J$ in $X_I$ whose
obstruction to compression is trivial, so Theorem A implies
the first part of Theorem B.
\S4 concerns the proof of Theorems D and E, which 
are  a consequence of Theorem B and 
 Milgram's EHP sequence.
In \S5 we prove Theorem A. The main tool in the proof is the relative
embedded thickening theorem of \cite{Klein3}. 
In \S6 we show that the embedding constructed
in \S3 has the correct collapse, thereby completing the proof of Theorem B.
In \S7 we prove the
stability of the obstruction $[\nu^+]$. 
In \S8 we prove Theorem F.

\subsection{Acknowledgments} This paper could not have been
written were it not for discussions I had with Tom Goodwillie
and Bill Richter. The proof of Theorem A was in part motivated
by  techniques employed  by Goodwillie to study the stability map in
relative pseudoisotopy theory. As I mentioned above, the first proof of
Theorem B is due to Richter. Also, the idea of the proof of \ref{embed-thick}
was aided by interaction with Richter. Thanks to 
 Andrew Ranicki for improvements in the exposition.
Lastly, I've benefited from  the papers of Bruce Williams.

\section{Preliminaries}
Our ground category is 
 $\Top$, the category  of compactly generated Hausdorff
spaces. This comes equipped
 with the structure of a Quillen  model category:
\begin{itemize}
\item The {\it  weak equivalences} are the  weak homotopy equivalences
(i.e., maps $X \to Y$ such that the associated
realization of its singular map $|S_\cdot X| \to |S_\cdot Y|$
is a homotopy equivalence). Weak equivalences are 
denoted $\overset\sim \to$.
 \item The {\it fibrations}, denoted $\twoheadrightarrow$,  are the Serre fibrations.
\item  The {\it cofibrations},
denoted $\rightarrowtail$,
are the `Serre cofibrations', i.e.,
inclusion maps given by a sequence of
cell attachments (i.e., relative 
cellular  inclusions) or retracts thereof.
\end{itemize}

Every object is fibrant.  
The cofibrant objects are the retracts of iterated cell attachments
built up from the empty space. Every object $Y$ comes equipped with
a functorial cofibrant approximation $Y^{\text{c}} 
\,\,{}^\sim\!\!\!\!\!\!\twoheadrightarrow Y$.

 A non-empty space is always  $(-1)$-connected. A connected
 space is $0$-connect\-ed, and is $r$-connected for some $r > 0$ if
 its homotopy groups vanish up through degree $r$, for any choice
 of basepoint. A map of non-empty spaces $X \to Y$ is called 
 {\it $r$-connected}
 if its homotopy fiber with respect to any choice of basepoint in $Y$
 is an $(r{-}1)$-connected space. 
 An $\infty$-connected map is a weak
 equivalence. 

A space is {\it homotopy finite} if it is homotopy equivalent to 
a finite CW complex.

A commutative square of cofibrant spaces
$$
\begin{CD}
A @>>> B\\
@VVV @VVV \\
C @>>> D
\end{CD}
$$
is {\it $r$-cocartesian} if the evident map
$C_0 \cup A_{[0,1]} \cup B_1 \to D$ (whose source
is a double mapping cylinder) is $r$-connected. More generally,
a square of not necessarily cofibrant spaces is 
$r$-cocartesian if it is after applying cofibrant approximation.
An $\infty$-cocartesian square is {\it cocartesian}. 
Dually,  the square is {\it $r$-cartesian} if the
map $A \to \holim{(B \to D \leftarrow C)}$ is $r$-connected.
An $\infty$-cartesian square is {\it cartesian}.

We introduce one last non-standard notation:
given a  map  of spaces $A \to B$, if no confusion arises 
we will often let $(\bar B,A)$ denote the pair
given by the mapping cylinder $B_0 \cup A_I$ with the inclusion
of $A_1$.

\subsection{Fiberwise spaces} \label{fib-space} For $X\in \Top$ an object,
  $R(X)$ will denote the category of retractive
 spaces, as in the introduction 
 (in another notation, not to be used here, $X\backslash {\Top}/X$).
 We will assume in what follows that $X$ is a cofibrant object of
 $\Top$.
 
 According to Quillen \cite{Quillen}, $R(X)$ inherits a model category
 structure arising from the one on $\Top$. 
Weak equivalences and fibrations are defined using
the forgetful functor $R(X) \to \Top$.
Cofibrations
 are those maps satisfying  the left lifting property with respect
 to the acyclic fibrations (the word `acyclic' is synonymous with
weak equivalence).
 
Any object  $Y \in R(X)$ comes 
equipped with a functorial cofibrant approximation
$Y^{\text c}\,\,\,{}^\sim\!\!\!\!\!\!\twoheadrightarrow  Y$
 and similarly, a functorial fibrant approximation
$Y\,\,\, {}^\sim\!\!\!\!\!\!\!\rightarrowtail
 Y^{\text{f}}$.

Given an object $Y \in R(X)$, define its {\it fiberwise suspension}
$\vs_X Y$ to be the object whose underlying space is obtained
by collapsing the subspace $X_I \subset \vs_X Y$ to $X$ (via the first
factor projection)
in the double mapping cylinder 
$X_0 \cup Y_I \cup X_1$. If $Y$ is cofibrant, then so
is its fiberwise suspension. We use the notation $\vs_X^j Y$ to denote
the $j$-fold iterated application of $\vs_X$ to $Y$.

The {\it homotopy category} of $R(X)$, denoted $\text{ho}R(X)$,
is the category whose objects are those of $R(X)$ and in which
the hom-set from an object $Y$ to an object $Z$ is given by
homotopy classes of morphisms $Y^{\text c} \to Z^{\text f}$.
This is denoted $[Y,Z]_X$; it is a pointed set.
The corresponding {\it stable} hom-set
is $\{Y,Z\}_X := \lim_j [\vs^j_X Y ,\vs_X^j Z]_X$. 

Obstruction theory in $\Top$ gives rise to an obstruction theory in $R(X)$.
Let $Z \in R(X)$ be an object. A commutative
diagram  
$$
\begin{CD} 
S^{j-1} @>>> Z \\
@VVV @VVV \\
D^j @>>> X
\end{CD}
$$
defines another object $Z \cup D^j_X$, whose underlying space
is $Z \cup_{S^{j-1}} D^j$. This operation is called {\it attaching
a $j$-cell} to $Z$.

\begin{defn} An object $P \in R(X)$ 
has {\it dimension} $\le s$ 
if its fibrant approximation admits
 a factorization $X \rightarrowtail
P' \overset\sim\to P^{\text{\rm f}}$ such that $P'$ is obtained from $X$ 
by attaching cells of dimension $\le s$.

A morphism $Y \to Z$ is {\it $r$-connected} if it
is  $r$-connected as a map of spaces. In particular,
an object $Y$ is {\it $r$-connected}
 if its structure map $X \to Y$
is.
\end{defn}

\begin{lem} \label{obs-theory} Let $Y \to Z$ be $r$-connected morphism
of $R(X)$ and suppose
that $P$ has dimension $\le r$. Then the induced map of homotopy
sets
$$
[P,Y]_X \to [P,Z]_X
$$
is surjective. It is also injective if $P$ has dimension $\le r{-}1$.
\end{lem}

This is essentially \cite[9.2]{James:Fiber-Hty}.

\subsection{The stable range}
The Freudenthal theorem measures the extent to which 
fiberwise suspension is an isomorphism on the level of fiberwise
homotopy classes.

\begin{lem} 
\label{freud}{\rm (James \cite[9.3]{James:Fiber-Hty}).}  If $Y,Z \in R(X)$ 
cofibrant 
objects such that $Z$ $r$-connected and  $Y$ has dimension $\le 2r {+} 1$, 
then fiberwise suspension gives a surjection of pointed sets
$$
[Y,Z]_X \to [\vs_X Y, \vs_X Z]_X \, .
$$
This surjection is an isomorphism whenever $Y$ has dimension $\le 2r$.
\end{lem}

\subsection{Poincar\'e spaces}  \label{pd-spaces}
In this paper, a {\it Poincar\'e space} $X$ of dimension $n$ is a 
pair $(X,\partial X)$ such that $X$ and $\partial X$ are cofibrant and
homotopy finite, $\partial X \to X$ is a cofibration, and $X$ satisfies 
{\it Poincar\'e duality:}
\begin{itemize} 
\item  there exists a local system of abelian groups
 $\cal L$ of rank one 
defined on $X$, and a fundamental class
$[X] \in H_n(X,\partial X; \cal L)$ such that 
the cap product homomorphisms 
$$\cap [X]\:H^*(X;M) \to H_{n{-}*}(X,\partial X;{\cal L} \otimes M)
$$ 
and 
$$\cap [\partial X]\:
H^*(\partial X;N) \to H_{n{-}*{-}1}(\partial X;{\cal L}_{|\partial X} \otimes N)
$$
are isomorphisms, where 
$[\partial X] \in H_{n{-}1}(\partial X;{\cal L}_{|\partial X})$ is the image of
$[X]$ under the connecting homomorphism in the 
 homology exact sequence of  the pair $(X,\partial X)$, and
$M$ ($N$) is any local system  on $X$ (resp.\ on $\partial X$)
(compare \cite{Klein}, \cite{Wall3}).
\end{itemize}

If $(X,\partial X)$ is a pair such that $\partial X \to X$ is $2$-connected,
then the first duality isomorphism implies the second one 
(cf.\  \cite[2.1]{Klein}). In these circumstances,
$X$ is $n$-dimensional Poincar\'e if and only if 
$X_I$ is $(n{+}1)$-dimensional Poincar\'e.

\subsection{Embeddings}\label{defn-embed}
Let $M$ and $X$ a Poincar\'e spaces of dimension $n$, where $X$ is  connected.
An {\it embedding} of $M$ in $X$ is a commutative cocartesian square
of cofibrant homotopy finite spaces
$$
\begin{CD}
\partial M @>>> C \\
@V \text{incl.} VV @VVg V \\
M @>>f> X
\end{CD}
$$
together with a factorization of the inclusion
$\partial X \to C \to X$, such that
$(M,\partial M)$ and $(\bar C,\partial M \amalg \partial X)$
satisfy Poincar\'e duality with respect to the fundamental
classes 
obtained by taking the image of a fundamental class for $X$ 
under the homomorphisms
$$
H_n(X,\partial X;{\cal L}) \to 
H_n(\bar X,C;{\cal L}) \cong 
H_n(M,\partial M;f^*{\cal L})
$$
and 
$$
H_n(X,\partial X;{\cal L})
\to
H_n(X, M \amalg \partial X;{\cal L}) 
\cong 
H_n(\bar C,\partial M \amalg \partial X; g^*{\cal L})\, .
$$
If $\codim M \ge 3$ then one only need verify the
compatibility of  fundamental classes for $M$ (see
\cite[2.3]{Klein}).

The space $C$ is called the {\it complement}, and 
$f\: M \to X$ is the {\it underlying map} of the embedding.

The {\it decompression} of an embedding of  $M$ in $X$ is the
embedding of  $M_I$ in $X_I$ defined by the diagram
$$
\begin{CD}
\partial M_I  @>>> W \\
@VVV @VVV \\
M_I @>>> X_I
\end{CD}
$$
where $W = X_0 \cup C_I \cup X_1$ is ({\it unreduced}) fiberwise suspension,
and the factorization $\partial X_I \to W \to  X_I$ 
is  evident.

Two embeddings from $M$ to $X$ with complements
$C_0$ and $C_1$  are {\it elementary concordant} if there exists
a diagram of pairs
$$
\begin{diagram}
\node{(\partial M_I,\partial M_0 \amalg \partial M_1)}
 \arrow{e} \arrow{s} 
\node{(W,C_0 \cup (\partial X)_I \cup  C_1)}\arrow{s}\\
\node{(M_I, M_0 \amalg M_1)}\arrow{e} \node{(X_I, \partial X_I)}
\end{diagram}
$$
in which
 each associated  diagram of spaces
$$
\begin{diagram}
\node{\partial M_I}\arrow{e} \arrow{s} \node{W}\arrow{s}\\
\node{M_I}\arrow{e} \node{X_I}
\end{diagram}
\qquad \text{and} \qquad 
\begin{diagram}
\node{\partial M_0 \amalg \partial M_1}
  \arrow{e}\arrow{s} \node{C_0 \cup (\partial X)_I \cup C_1} \arrow{s}\\
\node{M_0 \amalg M_1}  \arrow{e} \node{\partial X_I}
\end{diagram}
$$
is cocartesian  (the latter  of these is obtained
from the disjoint union of the embedding diagrams 
using the inclusion $\partial X_0 \amalg \partial X_1 \subset \partial X_I$).
Moreover, the maps $C_i \to W$ are required to be weak equivalences.
More generally, {\it concordance} is the equivalence relation generated by elementary concordance.

\subsection{Embedded thickenings}
Suppose that $K$ is a cofibrant
space which is homotopy equivalent to a finite connected CW complex
of dimension $\le k$. Let $f\: K \to X$ be a map, where 
$X^n$ is a connected Poincar\'e space of dimension $n$. A 
cocartesian square
$$
\begin{CD}
A @>>> C \\
@VVV @VVV \\
K @>f>> X
\end{CD}
$$
(in which $A$ and $C$ are cofibrant and  homotopy finite),
together with a factorization $\partial X \to C \to X$
is called an {\it embedded thickening} of $f$ if
\begin{itemize}
\item $(\bar K, A)$ gives an $n$-dimensional Poincar\'e space such that
$\codim \bar K \ge n{-}k$, and 
\item Replacing $K$ by $\bar K$ in the diagram yields an embedding
in the sense of \ref{defn-embed}.
\end{itemize}
An embedded thickening is what was
called a {\it PD embedding}
in the terminology of \cite{Klein}. In order to avoid
confusion, we have changed the name to
distinguish between the embeddings appearing in
this paper (where the boundary data are {\it a priori} given)
and the ones of \cite{Klein} (embeddings of complexes in Poincar\'e spaces).

\section{Proof of Theorem B (first part)}
We show how Theorem A can be used to construct an
embedding of $M$ in $X$ from an unstable fiberwise normal invariant.

Let $\alpha_f \in [X\dbslash \partial X,M\dbslash \partial M]_X$ be
an unstable fiberwise normal invariant associated to an immersion
$f\:M \to X$.
Based on a construction
of Browder \cite{Browder3}  we will associate a Poincar\'e embedding
of $M_J$ in $X_I$. 

For this section only, let us agree that $M\dbslash \partial M$ now means
the object of $R(X)$ whose underlying space
is $X_0 \cup (\partial M)_I \cup M_1$ (the formulation provided
in the introduction differs from this description by a canonical weak
equivalence). Similarly, let $X\dbslash \partial X$ now mean
$X_0 \cup (\partial X)_I \cup X_1$. Let $h\:J \to I$ be the
homeomorphism $t \mapsto 3t {-} 1$. 

Then there is a commutative diagram of spaces
$$
\begin{diagram}
\node{\partial M_J}\arrow{s} \arrow{e}\node{M\dbslash \partial M} \arrow{s}
 \\
\node{M_J} \arrow{e} \node{X_I}
\end{diagram}
$$
in which the top  arrow is defined by  
  $$
\begin{CD}
M_{1/3} \cup (\partial M)_J \cup M_{2/3} @>\id {\times} h >> 
M_{0} \cup (\partial M)_I \cup M_1 @>f \cup {\id} \cup \id  >> 
X_0 \cup (\partial M)_I \cup M_1 \, ,
\end{CD}
$$
the bottom arrow is $f {\times} h$, and the vertical arrows
are evident. This diagram is cocartesian. 
In what follows, we must replace $M\dbslash \partial M$
in the diagram with its fibrant approximation
 $(M\dbslash \partial M)^{\text{f}}$. 
Assume that this has been done.

The Poincar\'e
 boundary for $X_I$ is $X\dbslash \partial X$; it
 factors through $(M\dbslash \partial M)^{\text{f}}$ via a representative for
 $\alpha_f$. This 
 defines the embedding of $M_J$ 
 in $X_I$. In particular, the complement of this embedding is 
  $(M\dbslash \partial M)^{\text{f}}$.
 
 Applying Theorem A, we see that the given embedding
 compresses to a embedding of $M$ in $X$
if and only if $[\nu^+] \in [M^+,M\dbslash \partial M]_X$
 is the trivial element.
 But by construction, $[\nu^+]$ is the fiberwise homotopy
 class determined  by making the composite (fiberwise) map
 $$
 M_{1/3}   \to M_{1/3} \cup (\partial M)_J \cup M_{2/3} \to 
X_0 \cup (\partial M)_I \cup  M_1
 $$
 ``based'' (i.e.,  add on a disjoint copy of $X$ to $M_{1/3}$).
 The composite clearly factors through the ``basepoint'' $X_0
 \subset X_0 \cup (\partial M)_I \cup  M_1$, 
so $[\nu^+]$ is the trivial element.

It remains to check that the collapse of the embedding
of $M$ in $X$ equals $\alpha_f$. This is not a formal consequence
of Theorem A, but rather, a consequence of the construction of the
particular embedding in the proof of Theorem A
contained in \S5 below. For this reason, we defer the proof of this
until  \S6.
 
 \section{Proof of Theorems D and E}

\begin{proof}[Proof of Theorem D] 
We first explain the idea of the proof while ignoring technical details.
There is a commutative diagram of $R(X)$
$$
\begin{diagram}
\node{M\dbslash \partial M}\arrow{s} \arrow{e}\node{(M/\partial M){\times} X}\arrow{s} \\
\node{Q_X (M\dbslash \partial M)} \arrow{e} 
\node{Q (M/\partial M)  {\times} X} 
\end{diagram}
$$
in which
\begin{itemize}
\item  $(M/\partial M) {\times} X$ has structure maps
given by the second factor projection and the inclusion $* {\times} X
\subset (M/\partial M) {\times} X$. 
\item The morphism $M\dbslash \partial M \to (M/\partial M){\times} X$
 is given by the quotient map $M\dbslash \partial M \to  M/\partial M$
together with the retraction
$M\dbslash \partial M \to X$. 
\item $Q_X$ means the fiberwise version of stable homotopy, and the
bottom map of the diagram is defined in a way similar to the top map.
\item The vertical maps are defined by means of the natural transformation
from the identity to (fiberwise) stable homotopy.
\end{itemize}
\
Ignoring for the moment the issue of homotopy invariance of the terms
in the diagram, it will follow by an argument  sketched below
that the  diagram is $n$-cartesian.
Assuming this the argument proceeds as follows:

The fiberwise stable homotopy class $\alpha^{\text{st}}$
is represented by a morphism $X\dbslash \partial X
\to Q_X (M\dbslash \partial M)$ and the  homotopy
class $\beta$ is represented by a morphism $X\dbslash \partial X
\to (M/\partial M) {\times} X$. 
Up to fiberwise
homotopy the maps are compatible with
the diagram.
By \ref{obs-theory} applied to 
the $n$-connected morphism  
$$
M\dbslash \partial M \,\, \to \,\,
\holim (Q_X (M\dbslash \partial M) \to Q(M/\partial M){\times} X
\leftarrow (M/\partial M){\times} X)
$$
there is an unstable fiberwise normal invariant
$\alpha \in [X\dbslash \partial X, M\dbslash \partial M]_X$.
Theorem D now follows by application of Theorem B.
\medskip

We now proceed to establish the degree to which the square
is cartesian. First of all, we replace the square
by an equivalent one which is homotopy invariant  
(for the extent to which $Q_X$ is a homotopy invariant functor
is unclear, even for objects which are fibrant
and cofibrant).

Choose a basepoint for $X$. 
Since $X$ is a connected cofibrant space, there is a homotopy
equivalence $X \simeq BG$ where $G$ is the geometric realization of
the  simplicial set given which is the Kan loop group of
the total singular complex of $X$. Here, we think of $G$
as a topological group object in $\Top$.
In what follows, we will assume $X$ is $BG$.

Let $R^G(*)$ denote the category of based $G$-spaces.
This admits the structure of a model category in
which a morphism $Y \to Z$ is a weak equivalence if (and only if)
it is a weak homotopy equivalence of spaces. 
Every object is fibrant and the cofibrant objects
are the retracts of free based $G$-CW complexes. In fact,
the homotopy categories of $R^G(*)$ and
$R(BG)$ are equivalent (but we will not require this.)

Let $M^{\sim}$ denote the pullback of $M \to BG \leftarrow EG$.
Then $M^\sim /\partial M^\sim $ is an
object of $R^G(*)$. We recover
$M\dbslash \partial M \in R(BG)$ up to weak equivalence 
by taking the Borel construction 
$(M^\sim /\partial M^\sim) {\times}_G EG$. We recover $M/\partial M$
as the {\it homotopy orbits} (i.e., reduced Borel construction)
$(M^\sim /\partial M^\sim)_{hG} := 
(M^\sim /\partial M^\sim) \smsh_G EG_+$. 
In its homotopy invariant formulation,  the square is now given by
the diagram of morphisms of $R(BG)$
\begin{equation}\label{n-cart-square}
\begin{diagram}
\node{(M^\sim /\partial M^\sim)^{\text{c}} {\times}_G EG}
 \arrow{e}\arrow{s} \node{
(M^\sim /\partial M^\sim)^{\text{c}}_{hG} {\times} BG} \arrow{s}\\
\node{Q ((M^\sim /\partial M^\sim)^{\text{c}})  {\times}_G EG}
\arrow{e} \node{Q ((M^\sim /\partial M^\sim)^{\text{c}}_{hG}) {\times} BG}
\end{diagram} 
\end{equation}
(here, for an object $Y \in R^G(*)$,
the object $Y^{\text{c}}$ denotes its cofibrant approximation).

Finally, we calculate the degree to which the square is cartesian. 
In what follows, set
 $N := M^\sim /\partial M^\sim$, and
note  that $N$ is $(n{-}p{-}1)$-connected.
The homotopy fiber of the left vertical map is the same thing
as the homotopy fiber of the map $N \to
Q N$. Denote this fiber by $F_1$. By  Milgram's EHP-sequence 
\cite[1.11]{Milgram}, 
there is
a $(3n{-}3p{-}3)$-connected
map $\vo ( N \smsh N )_{h{\Bbb Z}/2} \to F_1$. On the other hand
the homotopy fiber of the right vertical map is the same as the homotopy
fiber of the map $N_{hG} \to Q(N_{hG})$. If we denote this homotopy
fiber by $F_2$, it again follows by Milgram's EHP-sequence that there
is a $(3n{-}3p{-}3)$-connected
map $\vo (N_{hG} \smsh N_{hG})_{h{\Bbb Z}/2} \to F_2$. Moreover,
the square
$$
\begin{diagram}
\node{\vo ( N \smsh N )_{h{\Bbb Z}/2}} \arrow{e}\arrow{s}
 \node{\vo (N_{hG} \smsh N_{hG})_{h{\Bbb Z}/2}} \arrow{s}\\
\node{F_1} \arrow{e} \node{F_2} 
\end{diagram}
$$
is commutative.
The top map of the latter square is induced by the evident map
$N \smsh N \to  (N\smsh N)_{hG {\times} G}$. This last map is
easily checked to be $(2n{-}2p{+}[p/2])$-connected. Assembling this
information, it follows that the map $F_1 \to F_2$ is
$\min(3n{-}3p{-}4,2n{-}2p{+}[p/2]{-}1)$-connected. By hypothesis,
$3p{+}4 \le n$, so this connectivity is at least $n$. Consequently,
the square \eqref{n-cart-square} is $n$-cartesian, as claimed.
\end{proof}

\begin{proof}[Proof of Theorem E]
The proof is similar to the proof of Theorem D (where here
$P$ plays the role of $M$). Therefore,
we will only sketch the argument and leave it to the reader to
fill in the details.

As above, there is a diagram
$$
\begin{diagram}
\node{\widetilde P\dbslash 
\partial \widetilde P} \arrow{e}\arrow{s} 
\node{(\widetilde P/\partial \widetilde  P) 
{\times} \tilde X} \arrow{s} \\
\node{
Q_{\widetilde X} \widetilde P\dbslash 
\partial \widetilde P} \arrow{e} 
\node{Q(\widetilde P/\partial \widetilde P) 
{\times} \tilde X} 
\end{diagram}
$$
which one checks (by essentially the same argument) to be 
$(2p)$-cartesian. The fiberwise stable normal invariant
can be lifted to a fiberwise equivariant map
$ \widetilde X/\partial \widetilde X \to
   Q_{\widetilde X} \widetilde P\dbslash 
\partial \widetilde P$. The rest of the argument follows as in the
proof of Theorem B, substituting obstruction theory by equivariant
obstruction theory, and using the fact that the equivariant homotopy
dimension of $\widetilde X/\partial \widetilde X$ is $2p$.
\end{proof}

\section{Proof of Theorem A}

Our main  tool  will be the relative 
embedded thickening theorem of \cite{Klein3} (see also \cite{Klein} for the
absolute version). The statement of this result will require some
preparation.

Let $(K,L)$ be a cofibration pair in $\Top$. We assume for simplicity
that $K$ and $L$ are cofibrant spaces which are homotopy finite.
Write  $$
\dim(K,L) \le k
$$ 
if there exists a factorization $L \to K' \to K$
in which $K'$ is obtained from $L$ by attaching cells of dimension
$\le k$ and the map $K' \to K$ is a weak equivalence. 

Let 
$X$ be an $n$-dimensional Poincar\'e space. 

By a {\it relative embedded thickening} of $(K,L)$ in $(X,\partial X)$
we mean a commutative diagram of cofibration pairs
$$
\begin{CD}
(A_K,A_L) @>>> (C_K,C_L) \\
@VVV @VVV \\
(K,L) @>>> (X,\partial X)
\end{CD}
$$
having the following properties.
\begin{itemize}
\item Each space appearing in the diagram is cofibrant
and homotopy finite.
\item Each of the diagrams of spaces
$$
\begin{CD}
A_K @>>> C_K \\
@VVV @VVV \\
K @>>> X
\end{CD}
\qquad \text{and} \qquad 
\begin{CD}
A_L @>>> C_L \\
@VVV @VVV \\
L @>>> \partial X
\end{CD}
$$
is cocartesian and the latter of these diagrams is a embedded thickening
of $L$ in $\partial X$.
\item The image of the fundamental
class of $X$ with respect to the composite 
$$
H_n(X,\partial X)
\to  
H_n(\bar X, \partial X \cup_{C_L}  C_K ) \cong 
H_n(\bar K, L \cup_{A_L} A_K)
$$
gives $(\bar K, L \cup_{A_L}A_K)$ the structure of 
an $n$-dimensional Poincar\'e space (here, coefficients are given
by pulling back the orientation bundle for $(X,\partial X)$).
Similarly, $(\bar C_K, C_L \cup_{A_L} A_K)$ has the structure
of a Poincar\'e space with  fundamental class induced from $X$.
\item The map $A_K \to K$ is $(n{-}k{-}1)$-connected.
\end{itemize}
The decomposition of $(X,\partial X)$  
is depicted in figure 1 below.
\medskip

\centerline{
\epsfxsize=2.5in
\epsfbox{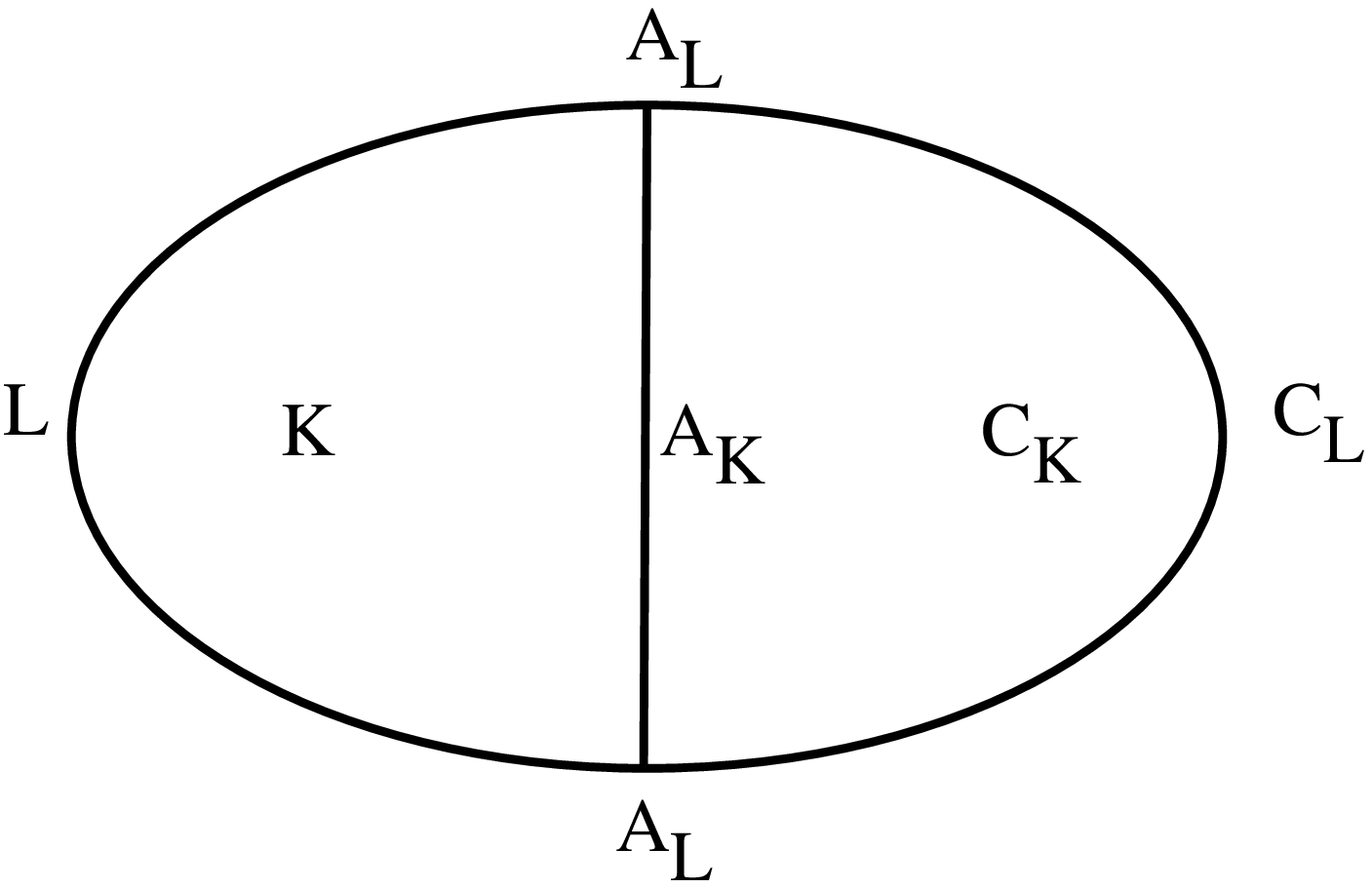}}
%\vspace{1.3in}
\centerline{\small Figure 1}
\medskip

Now let $f\:(K,L) \to (X,\partial X)$ be a map with
$\dim(K,L) \le k$ and suppose
that the restriction $f_{|L}\:L \to \partial  X$ 
embedded thickens.
The main theorem of \cite{Klein3} is 

\begin{thm} \label{rel-pd}
Assume  $k\le n{-}3$ and  $f\:K \to X$ is $(2k{-}n{+}2)$-connected.
Then there exists a relative embedded thickening
 of $f\:(K,L) \to (X,\partial X)$
extending  the given embedded thickening of $f_{|L}\:L \to \partial X$.
\end{thm}

\begin{rem} The above  is the
Poincar\'e version of the relative embedded thickening theorem of 
Hodgson \cite{Hodgson2}, with a loss of one dimension.
\end{rem}

We now begin the proof of Theorem A. Assume   
$[\nu^+]\in [M^+,W]_X$ is trivial, where $W$ is the complement of an
embedding of $M_J$ in 
$X_I$.  We may also 
assume without loss in generality that $W \in R(X)$ is fibrant.
A choice of  fiberwise null-homotopy may be thought of as a family of maps
$\nu_t\: M_t \to W$ for $t \in [0,1/3]$ which commute with projection to $X$ 
such that $\nu = \nu_{1/3}$ and $\nu_0$ factors through 
$X_0 \to W$.

This null-homotopy
gives rise to  a map of pairs
$$
(X_0 \cup M_{[0,1/3]}, X_0  \amalg M_{1/3})
\to
(W,\partial W)
$$
in which  $X_0 \cup M_{[0,1/3]}$ is the mapping cylinder
of the map $M_{1/3}  \to X$. These circumstances
 are depicted in figure 2.
\medskip

\centerline{
\epsfxsize=1.5in
\epsfbox{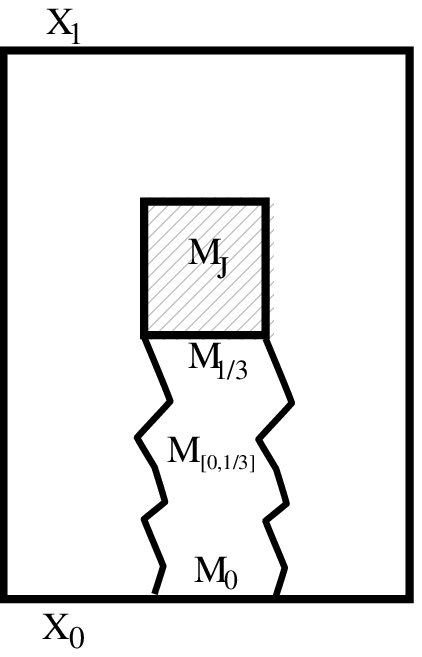}}
\centerline{\small Figure 2}
\medskip

The restricted map of spaces
$$
X \amalg M \to \partial W
$$ is already embedded thickened
(here, $\partial W = \partial X_I \amalg \partial M_J$).
This embedded thickening is given by the cocartesian square
$$
\begin{diagram}
\node{\partial X_0 \amalg \partial M_{1/3}}  \arrow{e}\arrow{s}   
\node{((\partial X)_I \cup  X_1 )
\amalg ((\partial M)_J \cup M_{2/3})}\arrow{s} \\
\node{X_0 \amalg M_{1/3}} \arrow{e} \node{\partial W\, .}
\end{diagram}
$$

The map
$$
X_0 \cup M_{[0,1/3]}  \to W
$$
is $(n{-}p{-}1)$-connected (since it, followed by the map 
$W \to X_I$ is a weak equivalence, and the 
latter map is $(n{-}p)$-connected).
Moreover, the pair $(X_0 \cup M_{[0,1/3]}, X  \amalg M)$ has relative dimension $\le p{+}1$.

Since $n{-}p{-}1 \ge 2(p{+}1) - (n {+}1) + 2$ if and only if
$2n \ge 3p {+} 4$, 
by \ref{rel-pd}   there exists a relative embedded thickening of 
$$
(X_0 \cup M_{[0,1/3]}, X  \amalg M) \to
(W,\partial W)
$$ 
which extends the given embedded thickening of 
$X \amalg M\to \partial W$. Thus we have a diagram of pairs (cf.\ fig.\ 3)
$$
\begin{diagram}
\node{(A,\partial X_0  \amalg \partial M_{1/3})} \arrow{e}\arrow{s} 
\node{(C,((\partial X)_I \cup  X_1 )
\amalg ((\partial M)_J \cup M_{2/3}))} \arrow{s}\\
\node{(X_0 \cup M_{[0,1/3]},X \amalg M)} \arrow{e} \node{(W,\partial W)\, .}
\end{diagram}
$$
\medskip

\centerline{
\epsfxsize=1.5in
\epsfbox{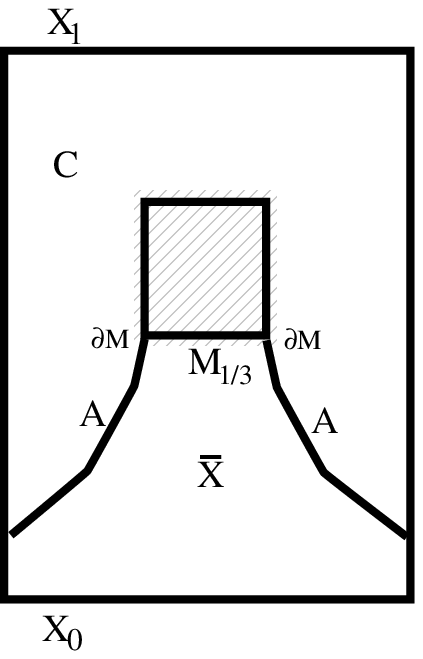}}
\vskip -.1in
\centerline{\small Figure 3}
\bigskip

Consider the associated commutative diagram
\begin{equation}\label{pd-diagram}
\begin{CD}
\partial M @>>> A \\
@VVV @VVV \\
M @>>> X
\end{CD}
\end{equation}
and note that there is an evident factorization of $\partial X \to X$
through the map $A \to X$.

To complete the proof of Theorem A, it suffices to show:

\begin{claim} The square  \eqref{pd-diagram} is an embedding of $M$ in
$X$. It induces the given embedding of $M_J$
in $X_I$ after decompression.
\end{claim}

To establish the claim, we first need to show that the square is cocartesian.
According to the definitions $X_0 \cup M_{[0,1/3]} \simeq X$ has an $n$-dimensional Poincar\'e boundary given by $X_0 \cup_{\partial X_0} (M_{1/3} \cup_{\partial M_{1/3}} A)$.  Application of Poincar\'e-Lefshetz duality
gives an isomorphism 
$$
H_*(\bar X,M \cup_{\partial M} A) \cong H^{n{+}1-*}(\bar X,X_0) = 0 \, .
$$
in all degrees, for any bundle of coefficients on $X$. Moreover, the 
map $M \cup_{\partial M} A \to \bar X$ induces an isomorphism
on fundamental groups (since $A \to X$
and $\partial M \to M$ are $2$-connected), so the square is cocartesian by 
application of Whitehead's theorem.

Secondly, a straightforward argument which we omit  shows that the
inclusion $X_1 \subset C$ is a weak equivalence. Consequently,
the composite $C \to W \to X$ is also a weak equivalence. Using this,
we have a chain of weak equivalences
$$
\vs_X A \overset\sim \to X \cup_A C \overset\sim \leftarrow W
$$
which is compatible with projection to $X_I$ and is relative
to $\partial X_I$. We infer that the decompression of \eqref{pd-diagram}
 yields the
embedding of $M_J$ in $M_I$ up to concordance. Compatibility of fundamental
classes is a consequence of the remarks
at the end of \ref{pd-spaces} and \ref{defn-embed}. 
This completes the proof of Theorem A.
\qed
\medskip

\section{Theorem B: completion of the proof}

Given a fiberwise unstable normal invariant 
$$
\alpha_f \in [ X\dbslash \partial X,M\dbslash \partial M]_X \, ,
$$
 we
constructed in \S3 an embedding of $M$ in $X$
by first associating an embedding of $M_J$ in $X_I$ and then
applying Theorem A (using the observation that the compression obstruction
of the latter embedding is trivial).

It remains to show that the collapse of this embedding 
coincides with $\alpha_f$.  We will give the argument
in the case when $\partial X = \emptyset$. The general
case, which is straightforward, will be left to the reader.

Returning to the proof of Theorem A and in particular fig.\ 3 above, note that the collapse  of the embedding
of $M$ in $X$ is the fiberwise homotopy class of the map 
$$
X_0 \amalg (A \cup_{\partial M} M) \to
\bar X \cup_{\partial M} M  \quad (\overset\sim \to  X \cup_{\partial M} M ) \, ,
$$
whose restriction to $X_0$ is given by the
inclusion $X_0 \to \bar X$ and the restriction to 
$A \cup_{\partial M} M$ is given by the amalgamation
of the map $A \to \bar X$ with the identity map of $M$.

Using fig.\ 3, we rewrite this as follows: 
consider the amalgamated union
$$
M' \quad := \quad 
(\partial M)_J  \cup_{\partial M_{2/3}} M_{2/3} \, .
$$
and write $X' := A \cup_{\partial M_{1/3}} M'$
(so $X'$ is identified with $X$ up to weak equivalence).
Then the
fiberwise homotopy class of the composite
$$
X_0 \amalg X' 
\to 
\bar X \cup_{\partial M_{1/3}} M' \overset \sim \to W
$$
represents the collapse of the embedding (recall that
$W$ is $M\dbslash \partial M$ made fibrant). 
Note there is an evident factorization
$X_0 \amalg X' \to X_0 \amalg C \to W$. 

On the other hand,  the composite
$$
X_0 \amalg X_1 \to X_0 \amalg C \to W
$$
induces $\alpha_f$. 

Consequently, the restrictions of the map
$X \amalg C \to W$ to $X_0 \amalg X'$ and
$X_0 \amalg X_1$ induce respectively
the collapse map of the embedding and $\alpha_f$. 

But the maps $X_1 \to C$ and $X' \to C$
are weak homotopy equivalences. Consequently, the map
$X_0 \amalg C \to W$ induces both the collapse of the embedding
of $M$ in $X$ and $\alpha_f$ on fiberwise homotopy.
Thus $\alpha_f$ coincides with the collapse.  
This completes the proof of Theorem B.\qed

\section{Stability of the obstruction}

To prove \ref{stab}, we apply  
\ref{freud} to the homotopy set $[M^+,W]_X$. 
Since $M$ is homotopy
equivalent to a complex of dimension $\le p$, we infer
that the object $M^+ \in R(X)$ has dimension $\le p$. On the other hand,
the connectivity of $W \in R(X)$ is one less than the connectivity
of the map $W \to X_I$, which in turn, is at least the connectivity
of the map $\partial M_J \to M_J$ since 
the former is the cobase change of the latter. 
But $\codim M_J \ge n{-} p{+}1 $, so $\partial M_J
\to M_J$ is $(n{-}p)$-connected. Hence $W \in R(X)$ is an 
$(n{-}p{-}1)$-connected
object. 

Consequently, \ref{freud}  implies that
$$
[M^+,W]_X \to [\vs_X M^+,\vs_X W]_X
$$
is an isomorphism whenever $p \le 2(n{-}p{-}1)$, or equivalently, whenever
$3p {+} 2 \le 2n$.
Thus, the obstruction to compression is stable in the range of Theorem A
(with two dimensions to spare).\qed

\section{Proof of Theorem F}

In this section we show how Theorem A implies a partial improvement
of the main result of \cite{Klein}. Let $K$ be a cofibrant space 
which is homotopy equivalent to a connected CW complex of dimension $\le k$.
Let $X$ be a connected $n$-dimensional Poincar\'e space.

The main result of \cite{Klein} is

\begin{thm} \label{klein} 
Assume that $f\:K \to X$ is $(2k {-} n{+}2)$-connected and
$k \le n{-}3$. Then there exists an embedded thickening of $f$.
\end{thm}

Now we have the statement of Theorem F, which is an
improvement of \ref{klein} in the metastable range:

\begin{thm} \label{embed-thick} Assume 
$f\: K \to X$ is $(2k{-}n{+}1)$-connected,
$k\le n{-}3$ and $3k{+}4 \le 2n$. Then there exists an embedded
thickening of $f$.
\end{thm}

\begin{proof} By \ref{klein}, there exists
an embedded thickening of the composite
$f_I \: K \to X = X_0  \overset\subset\to X_I$. Let this be denoted
$$
\begin{CD}
A' @>>> W \\
@VVV @VVV \\
K @>>f_I > X_I\, .
\end{CD}
$$
Without loss in generality, we may take $A' \to K$ to
be a fibration.
By straightforward application of
the Blakers-Massey theorem \cite[p.\ 309]{Good}, this square
is $k$-cartesian. Let $P$ denote
the homotopy pullback of the diagram given by deleting $A'$.
Then the evident map $A' \to P$ is $k$-connected.

The maps $\id\: K \to K$ and $K \overset f \to X = X_0 \subset W$
are compatible up to homotopy when followed by the given maps
to $X_I$. Consequently, there is an induced  map $K \to P$. As
$A' \to P$ is $k$-connected, we obtain a factorization 
$K \to A' \to P$. Since $A' \to K$ is a fibration, the homotopy
lifting property plus the factorization yield a section
$\zeta\:K \to A'$. By construction, the composite
\begin{equation} \label{3}
K^+ = K \amalg X  \overset{\zeta \amalg \id_{X}} \longrightarrow A' \amalg X_0 \to W
\end{equation}
is fiberwise null homotopic.

The map $\zeta\: K \to A'$ is $(n{-}k{-}1)$-connected. By 
\ref{klein}, it embedded thickens since $n{-}k{-}1 \ge 2k {-} n {+} 2$
is equivalent to $3k{+}3 \le 2n$.
Let 
$$
\begin{CD}
A @>>> C \\
@VVV @VVV \\
K @>>> A'
\end{CD}
$$
be such an embedded thickening. 
We claim that
the composite $C \to A' \to K$ is a weak equivalence. 
To see this, first note that $C \to K$ is
$2$-connected, since it is the composite of the $(n{-}k{-}1)$-connected
map $C \to A'$ with the $(n{-}k)$-connected map $A' \to K$.
Also, by Poincar\'e-Lefshetz duality, we infer that
$$
H_*(\bar K,C) \cong H^{n{+}1{-}*}(\bar K,K) = 0
$$
in all degrees. Consequently,
$C \to K$ is a weak equivalence by the Whitehead theorem.

Let $(M,\partial M)$ denote the pair $(\bar K,A)$. Then the argument
of the last paragraph implies that $(M_I,\partial M_I)$
coincides with $(\bar K,A')$ up to homotopy. Furthermore, with respect to this
homotopy equivalence, the inclusion $M_0 \subset \partial M_I$ corresponds
to $\zeta\: K \to A'$.

Assembling these data, we have an embedding of $M_I$ in $X_I$ whose
obstruction $[\nu^+]$ vanishes by \eqref{3}. Applying Theorem A 
yields an embedded thickening of $f\:K \to X$.
\end{proof}

\Addressesr

\end{document}